\documentclass[10pt]{article}
\usepackage[dvips]{graphicx}
\usepackage{amssymb}
\usepackage{amsmath}

\usepackage{caption}
\usepackage{hyperref}
\usepackage{arcs}
\usepackage {xcolor} 
\usepackage{subfigure}
\usepackage{authblk}
\usepackage{calrsfs}

\newtheorem{theorem}{Theorem}

\newtheorem{Remark}{Remark}

\usepackage{graphicx} % Required for inserting images

\newenvironment{AMS}{\small\bf 2020 AMS subject classification: }{} 

\title{Random sampling and polynomial-free interpolation by Generalized MultiQuadrics}

\author{A. Sommariva} 
\author{M. Vianello}
\affil{University of Padova, Italy}

\date{\today}

\begin{document}

\maketitle

\begin{abstract}
We prove that interpolation matrices for Generalized MultiQuadrics (GMQ) of order greater than one are almost surely nonsingular without polynomial addition, in any dimension and with any continuous random distribution of sampling points. We also include a new class of generalized MultiQuadrics recently proposed by Buhmann and Ortmann.
\end{abstract}

\vskip0.5cm
\noindent
\begin{AMS}
{\rm 65D05,65D12.}
\end{AMS}

\vskip1cm

\section{Introduction}
Looking at the relevant literature on RBF interpolation, the only results ensuring unisolvence of CPD (Conditionally Positive Definite) RBF {\em without polynomial addition} essentially date back to almost 40 years ago. These concern distance functions and CPD RBF of order 1 like MultiQuadrics \cite{M86}, and univariate cubic splines \cite{BS87} which are of order 3. Such a theoretical problem has been clearly stated in the popular monograph \cite{F07}: {\em ``There is no result
that states that interpolation with Thin-Plate Splines (or any other strictly conditionally positive definite function of order $m\geq 2$) without the addition of
an appropriate degree $m-1$ polynomial is well-posed''.} With the noteworthy exception of \cite{BS87} which has been in some sense overlooked (probably because the result is purely univariate), the problem remains substantailly open in the deterministic sampling framework.
On the other hand, the fact that for example interpolation by Thin-Plate Splines without polynomial addition is  practically always unisolvent, has been  substantially recognized in the application literature, cf. e.g. \cite{P22} with the references therein. 

Recently indeed, some advancement has been made in the random sampling framework. In particular, it has been proved that polynomial-free interpolation by Polyharmonic Splines of any order, in any dimension and with any continuous distribution of random sampling points is a.s. (almost surely) unisolvent. The proof is essentially based on a careful analysis of the determinant of the interpolation matrix, together with the fact that Polyharmonic Splines are {\em real analytic functions} up to their centers, where they present a singularity of some derivative (with additional difficulties in the case of Radial Powers of odd integer exponent); cf. \cite{BSV23,SV24}. The key underlying property is a famous basic result of measure theory, namely that {\em the zero set of a not identically zero real analytic function on an open connected domain has null Lebesgue measure}; cf. e.g. \cite{M20} with the references therein. It is worth quoting that a similar probabilistic  approach was also recently used in the framework of other (almost everywhere) analytic functions spaces \cite{DASV23,XN23}, 
for example polynomial spaces or spaces of RBF with fixed center, while the results quoted above concern the more difficult case of RBF whose centers coincide with the sampling points.

In the present note we extend the polynomial-free unisolvence result by random sampling to {\em Generalized MultiQuadrics} (GMQ), that correspond to the radial functions
\begin{equation} \label{GMQ}
\phi(r)=(1+(\varepsilon r)^{2k})^\beta\;,\;\;1<\beta\not\in \mathbb{N}\;,\;k\in \mathbb{N}\;.
\end{equation}
 The scale $\varepsilon >0$ represents the so-called \emph{shape parameter} associated with RBF and is important to control the interpolation problem conditioning, cf. e.g. \cite{F07,LS23}. 
For $k=1$, it is well-known that with $\beta<0$ (Generalized Inverse Multiquadrics) they are Positive Definite and with $\beta>0$ they are CPD of order $\lceil\beta\rceil$ (for $k=1$ and $\beta=1/2$ they are the classical Hardy's MultiQuadrics); cf. \cite{CS09,F07,H71}.
For $k>1$ these radial functions are a generalization of a new family recently proposed by Buhmann and Ortmann in the case $\beta=1/2$, which are shown to be not positive definite and where the parameter $k$ allows enhanced polynomial reproduction in the framework of quasi-interpolation; cf. \cite{OB24}. To our knowledge, such a family with its generalization appears here for the first time in the framework of scattered interpolation.

The main difficulty to prove polynomial-free unisolvence with respect to the case of TPS, is that $\phi$ is here an analytic function with {\em no real singularities}, so in the proof we have to embed the problem in the complex framework to exploit the presence of {\em complex singularities} (namely, branch points). 

\section{Polynomial-free random GMQ interpolation}

Our main result is the following: 

\begin{theorem} 
Let $\Omega$ be an open connected subset of $\mathbb{R}^d$, {\color{black}$d\geq 1$}, and 
$\{\mathbf{x}_i\}_{i\geq 1}$ be a randomly distributed sequence on $\Omega$ with respect to any given probability density $\sigma(\mathbf{x})$, i.e. a point sequence produced by sampling a sequence of absolutely continuous random variables $\{X_i\}_{i\geq 1}$ which are i.i.d. (independent and identically distributed) in $\Omega$ with density 
$\sigma\in L^ 1_+(\Omega)$. Moreover, let
\begin{equation} \label{Vn}
V_n=[\phi_j(\mathbf{x}_i)]\;,\;\;\phi_j(\mathbf{x})=\phi(\|\mathbf{x}-\mathbf{x}_j\|_2)\;,\;1\leq i,j\leq n\;,
\end{equation}
be the interpolation matrix with respect to Generalized MultiQuadrics (\ref{GMQ}) of order $\lceil\beta\rceil>1$. 
Then, for every $n\geq 1$ the points $\mathbf{x}_1,\dots,\mathbf{x}_n$ are a.s.(almost surely) distinct and the matrix $V_n$ is a.s. nonsingular. 
\end{theorem}
\vskip0.3cm
\noindent{\bf Proof.} The proof proceeds by complete induction on $n$. The induction base with $n=1$ is trivial, since $det(V_1)=\phi_1(\mathbf{x}_1)=\phi(0)=1$. 
Let us assume as inductive hypothesis that $det(V_k)$ is a.s. nonzero for 
$1\leq k\leq n$. 
First, $\mathbf{x}_{n+1}$ is a.s. distinct from $\mathbf{x}_1,\dots,\mathbf{x}_n$, since the probability that a random point falls in a finite set is clearly zero, as finite set having zero measure with respect to any distribution with density. 

Consider now the 
$(n+1)\times (n+1)$ matrix 
$$
A(\mathbf{x})=\left(\begin{array} {cccccc}
1 & \phi_2(\mathbf{x}_1)  & \cdots 
& \phi_{n-1}(\mathbf{x}_1) & \phi_{n}(\mathbf{x}_1)  & \phi_1(\mathbf{x})\\
\\
\phi_1(\mathbf{x}_2) & 1  & \cdots & \phi_{n-1}(\mathbf{x}_2) & \phi_{n}(\mathbf{x}_2)  & \phi_2(\mathbf{x})\\
\\
\vdots & \vdots  & \ddots & \vdots  & \vdots & \vdots\\
\\
\phi_1(\mathbf{x}_{n-1}) & \phi_2(\mathbf{x}_{n-1})  & \cdots & 1 & \phi_{n}(\mathbf{x}_{n-1})  & \phi_{n-1}(\mathbf{x})\\
\\
\phi_1(\mathbf{x}_{n}) & \phi_2(\mathbf{x}_{n})  & \cdots & \phi_{n-1}(\mathbf{x}_{n}) & 1  & \phi_{n}(\mathbf{x})\\
\\
\phi_1(\mathbf{x}) & \phi_2(\mathbf{x})  & \cdots & \phi_{n-1}(\mathbf{x}) & \phi_{n}(\mathbf{x}) & 1\\
\\
\end{array} \right)
$$
and observe that $A(\mathbf{x}_{n+1})=V_{n+1}$, 
since $\phi_j(\mathbf{x}_{n+1})=\phi_{n+1}(\mathbf{x}_{j})$.
Developing the determinant of $A(\mathbf{x})$ by Laplace rule along the last row, we get 
\begin{equation} \label{f}
f(\mathbf{x})=det(A(\mathbf{x}))=-det(V_{n-1})(\phi_n(\mathbf{x}))^2+a(\mathbf{x})\phi_n(\mathbf{x})+b(\mathbf{x})\;,
\end{equation}
where $a,b\in span\{\phi_1,\dots,\phi_{n-1}\}$. Notice that $f$ is an analytic function in $\mathbb{R}^d$, since it is defined by sum and products of functions in $span\{\phi_1,\dots,\phi_{n}\}$ and thus belongs to the function algebra generated by 
$\phi_1,\dots,\phi_{n}$; cf. \cite{KP02}. We claim that $f$ is a.s. not identically zero in $\Omega$. 

In fact, if $f$ were identically zero in $\Omega$, it would be identically zero also in $\mathbb{R}^d$, by the already quoted fundamental result that the zero set of a not identically zero real analytic function has null Lebesgue measure \cite{M20}, whereas $meas(\Omega)>0$. Then taking the line $\mathbf{x}(t)=\mathbf{x}_n+t\mathbf{u}$ where $\mathbf{u}=(u_1,\dots,u_d)$ is a given unit vector, we obtain that 
the real univariate function $f(\mathbf{x}(t))$ would be identically zero for $t\in \mathbb{R}$. Consequently, its analytic extension  to the complex plane, say $f(\mathbf{x}(z))$, would also be  identically zero for $z\in \mathbb{C}$.
Observe that $\|\mathbf{x}(z)-\mathbf{x}_j\|_2^2=\|\mathbf{x}_n+z\mathbf{u}-\mathbf{x}_j\|_2^2$ appearing in $\phi_j(\mathbf{x}(z))=\phi(\|\mathbf{x}(z)-\mathbf{x}_j\|_2)$, $1\leq j\leq n$, has to be seen not as the square of the complex 2-norm, but as the complex extension of the corresponding real function $$
\|\mathbf{x}_n+t\mathbf{u}-\mathbf{x}_j\|_2^2=\sum_{\ell=1}^d{((\mathbf{x}_n-\mathbf{x}_j)_\ell+tu_\ell)^2}
$$
$$
=\sum_{\ell=1}^d{[(\mathbf{x}_n-\mathbf{x}_j)_\ell^2+2tu_\ell(\mathbf{x}_n-\mathbf{x}_j)_\ell+t^2u_\ell^2]}
$$
$$
=t^2+2t\langle\mathbf{u},\mathbf{x}_n-\mathbf{x}_j\rangle+\|\mathbf{x}_n-\mathbf{x}_j\|_2^2\;,
$$
that is
$$
\|\mathbf{x}(z)-\mathbf{x}_j\|_2^2
=z^2+2z\langle\mathbf{u},\mathbf{x}_n-\mathbf{x}_j\rangle+\|\mathbf{x}_n-\mathbf{x}_j\|_2^2\;,
$$
where $\langle\cdot,\cdot\rangle$ denotes the Euclidean scalar product in $\mathbb{R}^d$.
Consequently
$$
\phi_j(\mathbf{x}(z))=\left(1+\varepsilon^{2k}\left(z^2+2z\langle\mathbf{u},\mathbf{x}_n-\mathbf{x}_j\rangle+\|\mathbf{x}_n-\mathbf{x}_j\|_2^2\right)^k\right)^\beta
$$
$$
=\exp\left(\beta\log\left(1+\varepsilon^{2k}\left(z^2+2z\langle\mathbf{u},\mathbf{x}_n-\mathbf{x}_j\rangle+\|\mathbf{x}_n-\mathbf{x}_j\|_2^2\right)^k\right)\right)\;,
$$
where $\log$ is the principal value of the complex logarithm. 

In particular
$$
\phi_n(\mathbf{x}(z))=\left(1+(\varepsilon z)^{2k}\right)^\beta
$$
has a branch point at {\color{black}$z_\ast=z_\ast(\varepsilon)=\varepsilon^{-1}e^{i\pi/(2k)}$, which is a root of the complex polynomial ($\varepsilon z)^{2k}+1$.} 
On the other hand, 
the functions $\phi_j(\mathbf{x}(z))$, $j<n$, are a.s. analytic at $z_\ast$, at least for a suitable choice of the unit vector $\mathbf{u}$.
To prove this fact, we have to distinguish the cases $k=1$ and $k>1$, {\color{black}and in the latter the cases $d\geq 2$ and $d=1$.} It is sufficient to check that the complex numbers 
$$
1+\varepsilon^{2k}\left(z_\ast^2+2z_\ast\langle\mathbf{u},\mathbf{x}_n-\mathbf{x}_j\rangle+\|\mathbf{x}_n-\mathbf{x}_j\|_2^2\right)^k\;,\;\;1\leq j<n\;,
$$
have positive real part, or even nonzero imaginary part, so that they cannot fall on the branch point or on the branch cut of the noninteger power $(\cdot)^\beta$. 

For $k=1$, it is sufficient to observe that for any $\mathbf{u}$ 
the complex numbers 
$$
1+\varepsilon^{2}\left((i/\varepsilon)^2+2i/\varepsilon\langle\mathbf{u},\mathbf{x}_n-\mathbf{x}_j\rangle+\|\mathbf{x}_n-\mathbf{x}_j\|_2^2\right)=\varepsilon^2\|\mathbf{x}_n-\mathbf{x}_j\|_2^2
+2i\varepsilon\langle\mathbf{u},\mathbf{x}_n-\mathbf{x}_j\rangle
$$
have a.s. positive real part since a.s. $\mathbf{x}_j\neq \mathbf{x}_n$ for $j<n$. 

For $k>1$, consider the complex-valued polynomial in $\mathbf{u}$ on the unit sphere
$$
p_j(\mathbf{u})=\varepsilon^{2k}\left(z_\ast^2+2z_\ast\langle\mathbf{u},\mathbf{x}_n-\mathbf{x}_j\rangle+\|\mathbf{x}_n-\mathbf{x}_j\|_2^2\right)^k\;.
$$

For {\color{black}$d\geq 2$}, taking $\mathbf{u}=\mathbf{u}_j$ orthogonal to $\mathbf{x}_n-\mathbf{x}_j$, we get $$p_j(\mathbf{u}_j)=\left((\varepsilon z_\ast)^2+\varepsilon^2\|\mathbf{x}_n-\mathbf{x}_j\|_2^2\right)^k=
\left(e^{i\pi/k}+\varepsilon^2\|\mathbf{x}_n-\mathbf{x}_j\|_2^2\right)^k\;,$$ which has nonzero imaginary part, since the (principal) argument of $e^{i\pi/k}+\varepsilon^2\|\mathbf{x}_n-\mathbf{x}_j\|_2^2$ is {\color{black}positive} and less than $\pi/k$ {\color{black}(this is geometrically seen at a glance by adding up the corresponding vectors in the complex plane)}. Then, also 
$p_j(\mathbf{u}_j)$ and thus $1+p_j(\mathbf{u}_j)$ have nonzero imaginary part, which means that the real polynomial $\mbox{Im}(1+p_j(\mathbf{u}))$ is not identically zero on the unit sphere and hence its zeros have null surface measure. The latter assertion can be proved for example by writing the unit vector $\mathbf{u}$ in spherical coordinates and observing that a nonzero real polynomial on the unit sphere becomes a multivariate trigonometric polynomial 
in the coordinate box, to which we can apply the fundamental result \cite{M20} on the zero set of real analytic functions.
Consequently, choosing a unit vector $\mathbf{u}$ not belonging to the union in $j$ of the zero sets of $\mbox{Im}(1 + p_j (u))$, the polynomial $1+p_j(\mathbf{u} )$ has nonzero imaginary part and thus the functions $\phi_j(\mathbf{x}(z))$ are analytic at $z_\ast$ for every $j<n$. 

{\color{black}For $d=1$, $\Omega$ is an open bounded interval and with no loss of generality we can interchange $x_n$ with $\max\{x_j:\,1\leq j\leq n\}$, which simply corresponds to interchange two rows and columns in the matrix $A(x)$. 
Choosing $u=1$ we get that the complex number $$p_j(1)=\left(e^{i\pi/k}+2\varepsilon (x_n-x_j)e^{i\pi/(2k)}+\varepsilon^2(x_n-x_j)^2\right)^k$$ has positive real part, since the (principal) argument of $e^{i\pi/k}+2\varepsilon (x_n-x_j)e^{i\pi/(2k)}+\varepsilon^2(x_n-x_j)^2$ is positive and less than $\pi/k$. Then also $1+p_j(1)$ has positive real part and again the functions $\phi_j(x(z))$ are analytic at $z_\ast$ for every $j<n$. 
}

Since for $k\geq 1$ the functions $\phi_j(\mathbf{x}(z))$ are a.s. analytic at $z_\ast$ for every $j<n$, also the functions $a(\mathbf{x}(z))$ and $b(\mathbf{x}(z))$ are both a.s. analytic at $z_\ast$, recalling that 
$a,b \in span\{\phi_1,\dots,\phi_{n-1}\}$.
Now, from $f(\mathbf{x}(z))\equiv 0$ we would get from (\ref{f})
$$
det(V_{n-1})\left(1+(\varepsilon z)^{2k}\right)^{2\beta}\equiv a(\mathbf{x}(z))\left(1+(\varepsilon z)^{2k}\right)^{\beta}+b(\mathbf{x}(z))\;.
$$

We show now that this identity leads to a contradiction for all real values of $\beta>1$, $\beta\not\in \mathbb{N}$. Let us begin with {\em rational} values of $\beta$, namely $\beta=p/q$ with $p$, $q$ relatively prime. If $q=2$, we would get
$$
det(V_{n-1})\left(1+(\varepsilon z)^{2k}\right)^{p}-b(\mathbf{x}(z))\equiv a(\mathbf{x}(z))\left(1+(\varepsilon z)^{2k}\right)^{p/2}
$$
which cannot hold, because the function on the left-hand side is analytic at $z_\ast$, whereas that on the right-hand side has a branch point of the complex square root there. If $q\neq 2$, let us restrict the $z$-variable to the segment $z=yz_\ast$ 
with $1-\delta<y<1$, where $\delta$ corresponds to a neighborhood  $|z-z_\ast|< \delta$ where $a(\mathbf{x}(z))$ and $b(\mathbf{x}(z))$ are analytic. Then we would have
$$
det(V_{n-1})\left(1-y^{2k}\right)^{2p/q}\equiv a(\mathbf{x}(yz_\ast))\left(1-y^{2k}\right)^{p/q}+b(\mathbf{x}(yz_\ast))\;,
$$
where we can write $1-y^{2k}=(1-y)(1+y+\dots+y^{2k-1})$. Moreover, by analiticilty and using the Taylor expansions of $a(\mathbf{x}(z))$ and $b(\mathbf{x}(z))$ at $z=z_\ast$ we would get $a(\mathbf{x}(yz_\ast))\sim c_s(1-y)^s$ and $b(\mathbf{x}(yz_\ast))\sim d_\nu(1-y)^\nu$ as $y\to 1^-$, 
where $c_s,d_\nu\neq 0$ and $s,\nu$ are the orders of the first respective nonvanishing derivatives.

Again, this cannot hold because if $\nu=0$ we would have distinct limits on the two sides as $y\to 1^-$. If $\nu>0$, looking at the infinitesimal orders for $y\to 1^-$, it would require that $2p/q=\min\{s+p/q,\nu\}$ i.e. either $p/q=s$ (that cannot hold 
since $\beta=p/q\not\in \mathbb{N}$) or $2p/q=\nu$ (that cannot hold since $p$ and $q$ are relatively prime and $q\neq 2$). Finally, consider the case of $\beta$ {\em irrational}. We would have
$$
det(V_{n-1})\left(1-y^{2k}\right)^{2\beta}\equiv a(\mathbf{x}(yz_\ast))\left(1-y^ {2k }\right)^{\beta}+b(\mathbf{x}(yz_\ast))
$$
and taking again the limits as 
$y\to 1^-$, this leads to a contradiction for $\nu=0$ (distinct limits on the two sides) and for $\nu>0$ because 
it would require that either $2\beta=s+\beta$ i.e. $\beta=s\in \mathbb{N}$, or $2\beta=\nu$ i.e. $\beta=\nu/2\in \mathbb{Q}$. 

We can now conclude the proof. If $f$ in (\ref{f}) is not identically zero, being analytic in $\Omega$ by \cite{M20}, its zero set, say $Z_f$, has null Lebesgue measure, and hence null measure with respect to any continuous probability distribution with density. On the other hand, $f(\mathbf{x}_{n+1})=det(V_{n+1})$. Then, we can write
$$
\mbox{prob}\{\mbox{det}(V_{n+1})=0\}=\mbox{prob}\{f(\mathbf{x}_{n+1})=0\}
$$
$$
=\mbox{prob}\{f\equiv 0\}
+\mbox{prob}\{f\not\equiv 0\;\&\;\mathbf{x}_{n+1}\in Z_f\}
=0+0=0\;,
$$
and the inductive step is completed.
\hspace{0.2cm} $\square$

\section*{Acknowledgements}

Work partially
supported by the
DOR funds 
of the University of Padova, and by the INdAM-GNCS 2024 Project ``Kernel and polynomial methods for approximation and integration: theory and application software".
This research has been accomplished within the RITA ``Research ITalian network on Approximation", the SIMAI Activity Group ANA\&A and the UMI Group TAA ``Approximation Theory and Applications".

\end{document}